\documentclass[graybox]{svmult}
\usepackage{mathptmx}
\usepackage{helvet}
\usepackage{courier}
\usepackage{makeidx}
\usepackage{graphicx}
\usepackage{multicol}
\usepackage[bottom]{footmisc}
\newcommand{\Teich}{\mbox{Teich}}
\newcommand{\Frob}{\mbox{Frob}}

\newcommand{\Q}{{\mathbb{Q}}}
\newcommand{\bbP}{{\mathbb{P}}}
\newcommand{\Z}{{\mathbb{Z}}}
\newcommand{\cH}{H}
\usepackage{amsmath,amssymb,amsfonts}
\newtheorem{conj}{Conjecture}
\makeindex
\begin{document}
\title*{Hypergeometric supercongruences}
\author{David P.\ Roberts and Fernando Rodriguez Villegas} 
\institute{ David P.\ Roberts \at University of Minnesota Morris, USA, 
\email{roberts@morris.umn.edu}
 \\ \\
Fernando Rodriguez Villegas \at The Abdus Salam International Centre for
Theoretical Physics, Italy,
\email{villegas@ictp.it}
}
\maketitle
\abstract{We discuss two related principles
for hypergeometric supercongrences, 
one related to accelerated convergence and 
the other to the vanishing of Hodge numbers. 
}
\section{Introduction}
\label{sec:1}

At the conference, we added two related principles to the study of
supercongruences involving the polynomials obtained by truncating
hypergeometric series.  By a {\it supercongruence} we mean a congruence
which somewhat unexpectedly remains valid when the prime modulus $p$ 
is replaced by $p^r$ for some integer $r>1$. We call $r$ the
{\it depth} of the supercongruence.

The first principle is that a supercongruence is the first instance of
a sequence of similar supercongruences, reflecting accelerated
convergence of certain Dwork quotients.  The second is that 
splittings of underlying motives can be viewed as the conceptual
source of supercongruences, with the depth of the congruence being
governed by the vanishing of Hodge numbers.

We present these principles here in a limited context, so that they
can be seen as clearly as possible. Let $\alpha=(\alpha_1$, \dots,
$\alpha_d)$ be a length $d$ vector of rational numbers in $(0,1)$ and
let $\beta = 1^d = (1,\dots,1)$.  We assume that that multiplication
by any integer coprime to the least common multiple $m$ of the
denominators of the $\alpha_i$'s preserves the multiset
$\{\alpha_1,\dots,\alpha_d\}$ modulo $\Z$.

The associated classical hypergeometric series and its $p$-power
truncations, for $p$ prime, are as follows.
\begin{align*}
F(\alpha,1^d|t) & := \sum_{k=0}^\infty \frac{(\alpha_1)_k \cdots
              (\alpha_d)_k}{k!^d} t^k,
& F_s(\alpha,1^d|t)&:= \sum_{k=0}^{p^s-1} \frac{(\alpha_1)_k \cdots
                (\alpha_d)_k}{k!^d} t^k.     
\end{align*}

Our starting point was the list $CY3$ of fourteen
$\alpha=(\alpha_1,\dots,\alpha_4)$ associated to certain families of
Calabi-Yau threefolds discussed in \cite{RV}.  Each has a
corresponding normalized Hecke eigenform $f=\sum a_n q^n$ of weight
four and trivial character.  For each, it was conjectured in \cite{RV}
that
\begin{equation}
\label{sc14}
F_1(\alpha,1^d|1) \equiv a_p \mbox{ mod } p^3, \qquad p\nmid ma_p.
\end{equation}
Some of these cases have been settled.  For example, the case
$\alpha = (1/5,2/5,3/5,4/5)$ was proved by McCarthy \cite{Mc}, the
corresponding modular form having level $25$~\cite{Sch}. Just before
submitting this note, Long, Tu, Yui, and Zudilin~\cite{LTYZ} announced
two different proofs of~\eqref{sc14} for all fourteen cases in $CY3$.

\section{Convergence to the unit root and Hodge gaps}
The two principles stem from observations about common behavior of the
examples in $CY3$.  The first observation is that each supercongruence
\eqref{sc14} seems to be part of a sequence.   Dwork proved~\cite{Dw}
that for $p\nmid m$
\begin{equation}
\label{dwork-congr}
\frac{F_{s+1}(\alpha,1^d\,|\,t)}{F_s(\alpha,1^d\,|\,t^p)} 
\equiv \frac{F_s(\alpha,1^d\,|\,t)}{F_{s-1}(\alpha,1^d\,|\,t^p)} \bmod
p^s, \qquad s\geq 0. 
\end{equation}
Moreover, the rational functions
$F_{s+1}(\alpha,1^d\,|\,t)/F_s(\alpha,1^d\,|\,t^p)$ converge as
$s\rightarrow \infty$ to a Krasner analytic function which can be
evaluated at a Teichm\"uller representative $\Teich(\tau)$ which is not
a zero of $F_1$ giving the {\it unit root}~$\gamma_p$ of the
corresponding local~$L$-series at~$p$.

For $\alpha\in CY3$, computations suggest
\begin{equation}
\label{super-congr-quintic-1}
\frac{F_s(\alpha,1^d\,|\,1)}{F_{s-1}(\alpha,1^d\,|\,1)}\equiv \gamma_p \bmod p^{3s}, \qquad p\nmid
ma_p, \qquad s> 0,
\end{equation}
where $\gamma_p\in \Z_p$ is the root of $T^2-a_pT+p^3$ not divisible by
$p$. Note that the case $s=1$ reduces to~\eqref{sc14}
since $\gamma_p\equiv a_p\bmod p^3$.

Our second observation is that the appearence of a congruence to a
power $p^{3s}$ as opposed to the expected $p^s$ is related to Hodge
theory.  Consider the hypergeometric family of motives
$\cH(\alpha,1^d\,|\, t)$ (see~\cite{MAGMA} for a computer
implementation).  For any $\tau \in \bbP^1(\Q)\setminus\{0,1,\infty\}$
the motive $\cH(\alpha,1^d\,|\, \tau)$ is defined over $\Q$, has
rank~$d$, weight $d-1$ and its only non-zero Hodge numbers are
$(h^{d-1,0},\ldots,h^{0,d-1})=(1,\ldots,1)$. When $\tau=1$  
there is a mild degeneration and the rank
drops to $d-1$.

For $\alpha\in CY3$, the motive for $\tau=1$ is the direct
sum, up to semi-simplification, of a Tate motive $\Q(-1)$ and the
motive $A=M(f)$ of the corresponding Hecke eigenform $f$ of weight
four. The Hodge numbers of $A$ are $(1,0,0,1)$.   
We view the gap of three between the initial $1$ and the next $1$ 
as explaining the supercongruences~\eqref{super-congr-quintic-1}.

\section{A congruence of depth five}
To illustrate our two observations further, we use the  
decomposition established in \cite[Cor.\ 2.1]{FOP} for the
case $\alpha=(1/2,1/2,1/2,1/2,1/2,1/2)$.  We learned at the
conference that this example was recently studied further by Osburn, Straub,
and Zudilin \cite{OSZ}, who proved~\eqref{super-congr} below for $s=1$
modulo $p^3$ and report that Mortenson conjectured it modulo $p^5$.

Again after semsimplifying, the motive $\cH(\alpha,1^6\,|\,1)$ has a
distinguished summand isomorphic to the Tate motive $\Q(-2)$ of
rank~$1$ and weight~$4$. The complement of this $\Q(-2)$ breaks up
into two pieces $A$ and $B$. They are both rank~$2$ motives of
weight~$5$.  Namely, $A=M(f_6)$ is the motive associated to the unique
normalized eigenform $f_6=\sum_{n\geq 1}a_n\, q^n$ of level $8$ and
weight $6$ and $B=M(f_4)(-1)$ is a Tate twist of the motive associated
to the unique normalized eigenform $f_4=\sum_{n\geq 1}b_n\, q^n$ of level $8$ and weight
$4$. The LMFDB~\cite{LMFDB} conveniently gives data on 
modular forms, including the $a_n$ and $b_n$ here.  

The trace of $\Frob_p$ on the full rank $5$ motive $\cH(\alpha,1^6\,|\,1)$
is given by
$$
a_p+b_pp+p^2.
$$
Numerically, we observe the following supercongruences 
\begin{equation}
\label{super-congr}
\frac{F_s(\alpha,1^d\,|\,1)}{F_{s-1}(\alpha,1^d\,|\,1)}\equiv \gamma_p\bmod p^{5s}, \qquad p\nmid
2a_p, \qquad s\geq 1,
\end{equation}
where $\gamma_p\in \Z_p$ is the root of $T^2-a_pT+p^5$ not divisible
by~$p$.

The Hodge numbers for $A$ and $B$ are $(1,0,0,0,0,1)$ and $(1,0,0,1)$
respectively, with the gap of five in the Hodge numbers for $A$ nicely
matching the exponent of the supercongruences.

\section{A summarizing conjecture}
We now state a conjecture that generalizes  the situations
discussed so far. 

\begin{conj}
\label{conj}
 For fixed $\tau=\pm1$, let $A$ be the unique submotive of
 $\cH(\alpha,1^d\,|\,\tau)$ 
  with $h^{0,d-1}(A)=1$ and let $r$ the smallest positive integer
  such that $h^{r,d-1-r}(A)=1$. 
For $p\nmid m$ such that $F_1(\alpha,1^d\,|\,\tau)\in \Z_p^\times$,
let $\gamma_p$ be the unit root of $A$. Then 
\begin{equation}
\label{conj-supercongr}
\frac{F_s(\alpha,1^d\,|\,\tau)}{F_{s-1}(\alpha,1^d\,|\,\tau)}
\equiv \gamma_p \bmod 
p^{rs}, \qquad s\geq 1.
\end{equation}
In particular, for $s=1$ we have
\begin{equation}
\label{conj-supercongr-1}
F_1(\alpha,1^d\,|\,\tau) \equiv a_p \bmod
p^r,
\end{equation}
where $a_p$ is the trace of $\mbox{\rm Frob}_p$ acting on $A$.
\end{conj} 

i) For generic $\alpha,\tau$ we expect $r=1$
and~\eqref{conj-supercongr} follows (see~\eqref{dwork-congr} and the
subsequent paragraph).  For the conjecture to predict $r>1$, the
motive has to split appropriately.

ii) For $\alpha = (1/2,\ldots,1/2)$ and $\tau = (-1)^d$ the motive
$\cH(\alpha,1^d\,|\,\tau)$ acquires an involution and we expect
$r=2$ for any $d \geq 7$; all numerical evidence is consistent with
this assertion.

iii) For large $d$ the unit roots involved are not in general related
to classical modular forms since the motives $A$ will  typically
have degrees greater than two.

\subsection*{Acknowledgments} 
We thank the organizers for the wonderful conference {\em
  Hypergeometric motives and Calabi-Yau differential equations}, held
at the Matrix Institute in Creswick, Australia, in January 2017.  DPR's research is supported by grant DMS-1601350
from the National Science Foundation. FRV would like to thank the AMSI
and the Australian Mathematical Society 
 for their financial support for his
participation in this conference. Special thanks go to our close
collaborator M.~Watkins for his continuous and significant
contributions to the hypergeometric motives project in general and for
the implementation of an associated package in {\it Magma}~\cite{MAGMA} in
particular.

\end{document}